\RequirePackage{fix-cm}
\documentclass[smallextended]{svjour3}       
\smartqed  
\usepackage{graphicx}
\usepackage{latexsym}
\usepackage{amssymb}
\usepackage{color}

\newcommand{\boldgreek}[1]{\mbox{\boldmath$#1$}}
\newcommand{\nfrac}[2]{{\textstyle \frac{#1}{#2} }}
\newcommand{\R}{I\kern-0.37emR}
\newcommand{\Q}{I\kern-0.37emP}
\newcommand{\E}{I\kern-0.37emE}
\newcommand{\ny}{n\rightarrow\infty}
\def\bgk{\boldgreek}
\newtheorem{thm}{Theorem}
\newtheorem{cor}{Corollary}
\newtheorem{rem}{Remark}

\def\bgk{\boldgreek}

\def\bbe{\bgk \beta}

 \journalname{Extremes ~
DOI 10.1007/s10687-015-0232-2\\}
\begin{document}

\title{Averaged extreme regression quantile\thanks{The research was supported by the grant GACR 15-00243S.}
}

\titlerunning{Averaged Extreme Regression Quantile}        

\author{Jana Jure\v{c}kov\'a        
}


\institute{J. Jure\v{c}kov\'a \at
              Charles University in Prague \\
              Faculty of Mathematic and Physics\\
              Tel.: +420-221913285\\
              Fax: +420-222323316\\
              \email{jurecko@karlin.mff.cuni.cz}           
}

\date{Received: date / Accepted: date}

\maketitle

\begin{abstract}
{Various events in the nature, economics and in other areas force us to combine the study of extremes with regression and other methods. A useful tool for reducing the role of nuisance regression, while we are interested in the shape or tails of the basic distribution, is provided by the \textit{averaged regression quantile} and namely by the  \textit{average extreme regression quantile}. Both are weighted means of regression quantile components, with weights depending on the regressors. 
      
Our primary interest is the averaged extreme regression quantile (AERQ), its structure, qualities and its applications, e.g. in investigation of a conditional loss given a value exogenous economic and market variables. AERQ has several interesting equivalent forms: While it is originally defined as an optimal solution of a specific linear programming problem, hence is a weighted mean of responses corresponding to the optimal base of the pertaining linear program, we give another equivalent form as a maximum residual of responses from a specific R-estimator of the slope components of regression parameter. The latter form shows that while AERQ equals to the maximum of some residuals of the responses, it has minimal possible perturbation by the regressors. 
 
Notice that these finite-sample results are true even for non-identically distributed model errors, e.g. under  heteroscedasticity. Moreover, the representations are formally true even when the errors are dependent - this all provokes a question of the right interpretation and of other possible  applications.      
}

\end{abstract}

\section{Introduction}
\label{intro}
\setcounter{equation}{0}
Consider the linear regression model
\begin{equation}
\label{1}
Y_{ni}=\beta_0+{\bf x}_{ni}^{\top}{\boldgreek\beta}+e_{ni}, \quad
i=1,\ldots,n
\end{equation}
with observations $Y_{n1},\ldots,Y_{nn},$
independent errors $e_{n1},\ldots,e_{nn},$ possibly 
non-identically distributed with unknown distribution functions $F_i, \; i=1,\ldots,n.$\\
The covariates ${\bf x}_{ni}=(x_{i1},\ldots,x_{ip})^{\top}, \; i=1,\ldots,n$ are random or nonrandom,
and
${\boldgreek\beta}^*=(\beta_0,\boldgreek\beta^{\top})^{\top}=(\beta_0,\beta_1,\ldots,\beta_p)^{\top}\in\mathbb R^{p+1}$ 
is an unknown parameter. We shall suppress the subscript $n$ whenever it does not
cause a confusion. Occasionally we use the notation $\mathbf x_{ni}^*=(1,x_{i1},\ldots,x_{ip})^{\top}, \; i=1,\ldots,n.$

The regression $\alpha$-quantile in model (\ref{1}) was originated by Koenker and Bassett \cite{KB1978}. 
The extreme regression quantile (extreme RQ)  was studied by Smith \cite{Smith1994} under \textit{\textit{\textit{iid}}} model errors 
(under the name ``nonregular regression''); he found  its
asymptotic distribution under heavy-tailed distribution $F$ and under
some conditions on the regressors. The extreme RQ under \textit{\textit{\textit{iid}}} model errors
was later studied  in \cite{Bassett1988}, \cite{Chern2005}, \cite{Jur2007a}, \cite{Jur2007b}, \cite{Knight1999}, \cite{Portnoy}, among others,
 under various
distributions and under miscellaneous regularity conditions. The asymptotic distributions of the intercept and slope components of the extreme RQ were 
derived in \cite{Jur2007a} and \cite{Jur2007b} under \textit{\textit{\textit{iid}}} errors for distributions with the tails lighter than $t$-distribution with 4 degrees of freedom and under some conditions on the regressors. 
In some cases the slope components of the extreme RQ even consistently estimate the slope parameters.

The regression $\alpha$-quantile $$\widehat{\boldgreek\beta}_n^*(\alpha)=\left(\hat{\beta}_{n0}(\alpha),(\widehat{\boldgreek\beta}_n(\alpha))^{\top}\right)^{\top}=
\left(\hat{\beta}_{n0}(\alpha),\hat{\beta}_{n1}(\alpha),\ldots,\hat{\beta}_{np}(\alpha)\right)^{\top}$$
is a $(p+1)$-dimensional vector defined by means of minimization
\begin{eqnarray}\label{RQ}
&&\widehat{\boldgreek\beta}_n^*(\alpha)=\arg\min_{\mathbf b\in{\R}^{p+1}}\Big\{\sum_{i=1}^n\Big[\alpha(Y_i-\mathbf x_i^{*\top}\mathbf b)^++(1-\alpha)(Y_i-\mathbf x_i^{*\top}\mathbf b)^-\Big]\Big\}\nonumber\\
&& \; \mbox{where } \; z^+=\max(z,0) \; \mbox{ and} \; z^-=\max(-z,0), \; z\in\mathbb R_1. 
\end{eqnarray} 
Particularly, if $\alpha=1,$ we get the extreme (maximal)
regression quantile $$\widehat{\boldgreek\beta}_n^*(1)=\arg\min_{\mathbf b\in{\R}^{p+1}}\sum_{i=1}^n(Y_i-\mathbf x_i^{*\top}\mathbf b)^+.$$ 
By (\ref{RQ}), $\widehat{\boldgreek\beta}_n^*(1)$  can be alternatively
described as any solution to the linear program:
\begin{equation} \label{maxdef}
\min_{\mathbf b^* \in {\R}^{p+1}} \, \sum_{i=1}^n \mathbf x_i^{*\top}\mathbf b^* \qquad %
{\mbox{s.t.}} \quad Y_i \leq \mathbf x_i^{*\top}\mathbf b^* , \quad %
i = 1,\ldots,n.
\end{equation}
If $\sum_{i=1}^nx_{ij}=0, \ j=1,\ldots,p,$ then we only minimize the intercept component $b_0,$
subject to the inequalities in (\ref{maxdef}). In the absence of regression, i.e. if\\ $\mathbf x_i=\mathbf 0, \; i=1,\ldots,n,$ the extreme regression quantile coincides with the extreme observation $Y_{n:n}.$ 

The weighted mean of components of $\widehat{\boldgreek\beta}_n^*(\alpha), \; 0\leq\alpha\leq 1$
\begin{equation}\label{22x}
B_n(\alpha)=\overline{\mathbf x}_n^{*\top}\widehat{\boldgreek\beta}_n^*(\alpha)=\widehat{\beta}_{n0}(\alpha)+\frac 1n\sum_{i=1}^n\sum_{j=1}^p x_{ij}\widehat{\beta}_j(\alpha)
\end{equation}
is called the averaged regression $\alpha$-quantile.  
The behavior of $B_n(\alpha)$ with $0<\alpha<1$ has been illustrated in \cite{Bassett1988} and \cite{KB1982}, and summarized in \cite{Koenker2005}. Interesting is the monotonicity of $B_n(\alpha)$ in $\alpha\in(0,1).$

 In the present paper, we are primarily interested in the structure, properties and possible use of averaged extreme regression quantile (AERQ) $B_n(1),$ corresponding to $\alpha=1.$ It is a useful tool when we are interested in the extreme behavior of our observations while they are affected by a regression.

The quantile regression was considered by a host of authors in investigation of a conditional loss given a value of a covariate $\mathbf x,$ involving exogenous economic and market variables (or otherwise the past observed
returns). Such problems we meet in the finance, but also in the insurance and in the social statistics. 

Particularly, in the context of the extreme regression quantile, we consider $Y_i$ as the loss variable 
under realization $\mathbf x_i$, representing economic and market (exogenous) variable. Assuming that the relation of the loss and covariates is described by model (\ref{1}), our primary goal is to find $\widehat{\bbe}^*=(\hat{\beta}_0, \widehat{\bbe})^{\top}$ leading to the minimum sum of positive residuals, those satisfying $Y_i>\mathbf x_i^{*\top}\bbe^*$. This reflects greater concern about underestimating losses Y,  
rather than overestimating, emphasized e.g. in \cite{RockMiranda}. 
Under continuous  distribution functions $F_i$ of $e_i, \; i=1,\ldots,n,$ with probability 1 we obtain  
\begin{eqnarray}\label{shortfall}
&&\overline{Y}_n-B_n(1)=\frac 1n \sum_{i=1}^n(Y_i-\mathbf x_i^{*\top}\widehat{\bbe}_n^*(1))\nonumber\\
&&=\frac 1n \sum_{i=1}^n(Y_i-\mathbf x_i^{*\top}\widehat{\bbe}_n^*(1))^+-\frac 1n \sum_{i=1}^n(Y_i-\mathbf x_i^{*\top}\widehat{\bbe}_n^*(1))^-\\
&&=-\frac 1n \sum_{i=1}^n(Y_i-\mathbf x_i^{*\top}\widehat{\bbe}_n^*(1))^-\leq 0.\nonumber
\end{eqnarray}
This can be interpreted as the expected conditional shortfall of the loss, which is non-positive under  $\bbe^*=\widehat{\bbe}_n^*(1).$ 
 We refer to \cite{Trindade}, \cite{CaiWang}, \cite{Rockafellar}, \cite{RockMiranda} for the discussions and various solutions, and to other papers cited in. For possible applications in the insurance we refer to \cite{Chan}.
 
 The structure of $B_n(1)$ is of interest.
Because $B_n(1)$ is an outcome of a linear programming, it is identical to a weighted average of $p+1$ among observations $(Y_1,\ldots,Y_n),$ 
corresponding to the optimal base of the linear program (\ref{maxdef}). 
We shall give another equivalent form of $B_n(1),$ which employs a suitable R-estimator of the slope components $\bbe.$ Combining two equivalent forms indicates that $B_n(1)$ never exceeds $\max_{i\leq n}Y_i,$ while relative to $\max_{i\leq n}Y_i$ it has minimum possible perturbations caused by the regressors, with minimum over values of $\bbe$.

\section{Finite-sample form of $B_n(1)$ based on optimal base} 
\label{sec:1}
\setcounter{equation}{0}
Consider the linear program (\ref{maxdef}); let $\{\mathbf x_{i_1},\ldots,\mathbf x_{i_{p+1}}\}$ be its optimal base and let $\{Y_{i_1},\ldots, Y_{i_{p+1}}\}$
be the corresponding responses. It follows from the linear programming theory  that $B_n(1)$  
equals to a weighted mean of $\{Y_{i_1},\ldots, Y_{i_{p+1}}\}$, with the weights based on the regressors. However, in this special case are all weights are positive, as it is specified in the following theorem.   
\begin{thm}\label{Theorem1}
Assume that the regression matrix $\mathbf X_n^*=\left[\begin{array}{l}  \mathbf x_{n1}^{*\top}\\  \ldots\\  \mathbf x_{nn}^{*\top}\end{array}\right]$ has full rank $p+1$ and that the distribution functions $F_1,\ldots,F_n$ of model errors are continuous and increasing in $(-\infty,\infty).$ Then with probability 1
\begin{eqnarray}\label{main0}
&&B_n(1)=\sum_{k=1}^{p+1} w_{k}Y_{i_k}, \quad w_k>0, \; k=1,\ldots,p+1, \; \sum_{k=1}^{p+1} w_{k}=1\\[2mm]
&&\overline{Y}_n\leq B_n(1)<\max_{i\leq n}Y_i\label{main01}
\end{eqnarray}
where $\overline{Y}_n=\frac 1n\sum_{i=1}^nY_i$ and the vector $\mathbf Y_n(1)=(Y_{i_1},\ldots,Y_{i_{p+1}})^{\top}$ 
corresponds to the optimal base of the linear program (\ref{maxdef}). 

The vector $\mathbf w=(w_1,\ldots,w_{p+1})^{\top}$ of coefficients in (\ref{main0}) has the form
\begin{equation}\label{w}
\mathbf w^{\top}=n^{-1}\mathbf 1_n^{\top}\mathbf X_n^*(\mathbf X_{n1}^*)^{-1};
\end{equation}
where $\mathbf X_{n1}^*$ is the submatrix of $\mathbf X_n^*$ with the rows $\mathbf x_{i_1}^{*\top},\ldots,\mathbf x_{i_{p+1}}^{*\top}.$
\end{thm}
\textbf{Proof.}
The dominance of $B_n(1)$ over the average $\overline{Y}_n$ follows directly from the definition of $B_n(1).$ The dominance of $\max_{i\leq n}Y_i$ over
$B_n(1)$ will follow from the fact proven below that the coefficients in (\ref{main0}) are all positive and sum up to 1.

To prove the identity in (\ref{main0}), consider the minimization (\ref{RQ}) for $0\leq\alpha\leq 1$ as a parametric linear programming. Then its dual program
is a parametric linear programming of the form
\begin{eqnarray}
\label{linprog}
&\mbox{ maximize \ } & {\bf Y}_{n}^{\top}\hat{{\bf a}}(\alpha)
\nonumber\\
&\mbox{under} & \mathbf X_n^{*\top}\hat{\mathbf a}(\alpha)=(1-\alpha)\mathbf X_n^{*\top}\mathbf 1_n^{\top}\\             
&& \hat{{\mathbf a}}(\alpha) \in [0,1]^n, \ 0\leq\alpha\leq 1.\nonumber
\end{eqnarray}
The components of the optimal solution $\hat{\bf a}(\alpha)=(\hat{a}_{n1}(\alpha),\ldots,\hat{a}_{nn}(\alpha))^{\top}$ of
(\ref{linprog}), called 
regression rank scores, were studied in \cite{GJ}, where it is shown that $\hat{a}_{ni}(\alpha)$ is a continuous, piecewise linear
function of $\alpha\in[0,1]$ and $\hat{a}_{ni}(0)=1, \; \hat{a}_{ni}(1)=0, \; i=1,\ldots,n.$ Moreover, $\hat{\bf a}(\alpha)$ is invariant in the sense that it does not change if $\mathbf Y$ is replaced with $\mathbf Y+\mathbf X_n\mathbf b, \; \forall \mathbf b^*\in\mathbb R^{p+1}$
(see \cite{GJ} for detail). 
The regression quantile $\widehat{\bbe}_n^*(\alpha)$ is a step function of $0\leq \alpha\leq 1.$ If $\alpha$ is a
continuity point of the regression quantile trajectory, then we have the following identities,  proven in \cite{JP2014}:
\begin{eqnarray}
\label{tail5}
&& B_n(\alpha)=\frac{1}{n}\sum_{i=1}^n{\bf x}_i^{*\top}\widehat{\boldgreek\beta}^*(\alpha)=
-\frac{1}{n}\sum_{i=1}^n Y_i\hat{a}_{ni}^{\prime}(\alpha)\\ 
\mbox{ or  }&&B_n(\alpha)-\overline{\mathbf x}_n^{*\top}\bbe^*=\frac{1}{n}\sum_{i=1}^n{\bf x}_i^{*\top}(\widehat{\boldgreek\beta}^*(\alpha)
-{\boldgreek\beta}^*)=
-\frac{1}{n}\sum_{i=1}^ne_i\hat{a}_{ni}^{\prime}(\alpha)\nonumber  \label{tail6}
\end{eqnarray} 
where $\hat{a}_{ni}^{\prime}(\alpha))=\frac{d}{d\alpha}\hat{a}_{ni}(\alpha).$ 
Moreover, (\ref{linprog}) implies
\begin{equation}
\label{tail2}
\sum_{i=1}^nx_{ij}^*\hat{a}_{ni}^{\prime}(\alpha) =-\sum_{i=1}^n x_{ij}^*, \ 1\leq i\leq n, \ 1\leq j \leq p, 
\end{equation}
and particularly,
\begin{equation}
\label{tail4a}
\frac 1n\sum_{i=1}^n \hat{a}_{ni}^{\prime}(\alpha)=-1. 
\end{equation}
Notice that 
$\hat{a}_{ni}^{\prime}(\alpha)\neq 0$  iff $\alpha$ is the point of continuity of $\widehat{\bbe}_n^*(\cdot)$ and $Y_i=\mathbf x_i^{*\top}\widehat{\boldgreek\beta}_n^*(\alpha), \;  i=1,\ldots,n.$ 
There are exactly $p+1$ such components pertaining to a fixed continuity point $\alpha;$ those for which  
$\mathbf x_i^*$ belongs to the optimal base of program (\ref{linprog}). The number of points of discontinuity of $\widehat\bbe_n(\cdot)$ is finite under a finite $n.$ Denote $1-\varepsilon_1<1$ the largest
 point 
such that, for $1-\varepsilon_1<\alpha<1$ all
$\hat{a}_{ni}(\alpha)$ are non-increasing, and let $1-\varepsilon_2\in(1-\varepsilon_1,1)$ be the largest point of discontinuity of $\widehat\bbe_n(\cdot).$ Then either $\hat{a}_{ni}^{\prime}(\alpha)=0$ or  $\frac{1}{\varepsilon_1}<-\hat{a}_{ni}^{\prime}(\alpha)<\frac{1}{\varepsilon_2}$ for $1-\varepsilon_1<\alpha<1;$ moreover,  $1-\varepsilon_1\leq 1-\frac 1n\leq 1-\varepsilon_2$ by (\ref{tail4a}).   
Hence, for $1-\varepsilon_1<\alpha<1$ are $\hat{a}_{ni}(\alpha)$ dec\-reas\-ing and $\hat{a}_{ni}^{\prime}(\alpha)< 0$ for exactly $p+1$ points $i_1,\ldots,i_{p+1},$ while $\hat{a}_{ni}^{\prime}(\alpha)= 0$ 
otherwise. 
 
Moreover, by 
(\ref{tail2})  and (\ref{tail4a})
\begin{eqnarray*}
\label{F}
&&\sum_{i=1}^n\hat{a}_{ni}^{\prime}(1)=-n\nonumber\\ 
&&\sum_{i=1}^n x_{ij}\hat{a}_{ni}^{\prime}(1)=-\sum_{i=1}^{n}x_{ij}, \
j=1,\ldots,p
\end{eqnarray*}
where $\hat{a}_{ni}^{\prime}(1)$ is the left-hand
derivative of $\hat{a}_{ni}(\alpha)$ at $\alpha=1.$
The identity 
(\ref{tail5}) particularly implies
\begin{equation}
\label{tail5a}
B_n(1)=-\frac{1}{n} \sum_{i=1}^n Y_i\hat{a}_{ni}^{\prime}(1),
\end{equation}
while exactly $p+1$ coefficients in (\ref{tail5a}) are different from zero, those corresponding to the base. Hence we have $B_n(1)=\sum_{k=1}^{p+1}w_kY_{i_k}.$

Notice that the inequalities among the constraints in (\ref{maxdef}) reconvert in equalities for just $p+1$ components of the optimal base $\mathbf x_{i_1}^*,\ldots,\mathbf x_{i_{p+1}}^*$.   
Let $\mathbf X_{n1}^*$ be the submatrix of $\mathbf X_n^*$ with the rows $\mathbf x_{i_1}^{*\top},\ldots,\mathbf x_{i_{p+1}}^{*\top}$ and let $(\hat{\mathbf a}^{\prime}(1))^{\top}=(\hat{a}_{i_1}^{\prime}(1),\ldots,\hat{a}_{i_{p+1}}^{\prime}(1)).$ Then $\mathbf X_{n1}^*$ is regular with probability 1 and  
$$\mathbf w^{\top}=(\hat{\mathbf a}^{\prime}(1))^{\top}=-\mathbf 1_n^{\top}\mathbf X_n^*(\mathbf X_{n1}^*)^{-1}.$$
This and (\ref{tail5a}) imply (\ref{main0}) and (\ref{w}). 
\hfill $\Box$
\begin{rem}
In the location case, where $Y_i=\beta_0+e_i, \; i=1,\ldots,n,$ we have $\mathbf X_n^*=\mathbf 1_n$ and the regression rank scores have the form (so called H\'ajek's rank scores, see \cite{Hajek65})
$$\hat{a}_{ni}(\alpha)=\left\{
\begin{array}{lll}
1 & \ldots & 0\leq \alpha\leq \frac{R_i-1}{n}\\[2mm]
R_i-n\alpha & \ldots & \frac{R_i-1}{n}<\alpha\leq \frac{R_i}{n}\\[2mm]
0 & \ldots & \frac{\R_i}{n}\leq\alpha\leq 1\\
\end{array}\right . , \qquad i=1,\ldots,n
$$
where $R_i$ is rank of $Y_i, \; i=1,\ldots,n$. Hence,
$$\hat{a}_{ni}^{\prime}(\alpha)=\left\{
\begin{array}{rll}
-n & \ldots & \frac{R_i-1}{n}<\alpha< \frac{R_i}{n}\\[2mm]
0  & \ldots & 0< \alpha<\frac{R_i-1}{n} \; \mbox{ or } \; \frac{\R_i}{n}<\alpha<1\\
\end{array}
\right .
$$ 
and particularly $\hat{a}_{ni}^{\prime}(1)=-n$ if $R_i=n,$ what corresponds to the maximal $Y_i,$ and $\hat{a}_{ni}^{\prime}(1)=0$ otherwise. Naturally, then $B_n(1)=\max_{i\leq n}Y_i.$ However, it follows from (\ref{main0}) that $B_n(1)<\max_{i\leq n}Y_i$ sharply in the presence of regression. The asymptotic behavior of $\max_{i\leq n}Y_{ni}-B_n(1)$ as $\ny$ in case of i.i.d. errors is of interest, and will be a subject of a further study.
\end{rem}

\section{Finite-sample expression of $B_n(1)$ using the R-estimator of slopes} 
\setcounter{equation}{0}
The identity in (\ref{main0}) does not {reveal any relation of $B_n(1)$} to the extreme $Y_{n:n}$ or to some extreme residuals, because the optimal base of the linear program can be found only numerically, not explicitly.
An alternative form of the extreme regression quantile, considered by the author  in \cite{Jur2007a}, treats the slope and intercept components in model (\ref{1}) se\-pa\-ra\-te\-ly: We start with a specific $R$-estimate $\widetilde{\bbe}_{nR}^+$ of the slope components $\bbe,$ based on the ranks of residuals $Y_i-\mathbf x_i^{\top}\mathbf b, \; i=1,\ldots,n$. 
This $R$-estimate $\widetilde{\bbe}_{nR}^+$ is defined as the minimizer of
the Jaeckel \cite{Jaeckel1972} measure of the rank dispersion
\begin{eqnarray}
\label{11da}
\widetilde{\bbe}_{nR}^+&=&\arg\min\{\mathcal D_n({\bf b}): \; \mathbf b\in{\R}^p\},\\[2mm]
{\mathcal D}_n({\bf b})&=&\sum_{i=1}^n(Y_i-{\bf x}_i^{\top}{\bf b})%
\varphi_{n}\left(\frac{R_{ni}(Y_i-{\mathbf x}_{i}^{\top}{\bf b})}{n+1}\right)\nonumber
\end{eqnarray}
with the ``extreme'' score function
\begin{equation}
\label{8b}
\varphi_{n}(u)=I[u\geq 1-\nfrac 1n]-\nfrac 1n, \quad 0\leq u\leq 1,
\end{equation}
where
$R_{ni}(Y_i-{\mathbf x}_{i}^{\top}{\bf b})$ is the rank of the residual
$Y_i-{\mathbf x}_{i}^{\top}{\bf b}$
among\\ $Y_1-\mathbf x_{1}^{\top}{\bf b},%
\ldots,Y_n-\mathbf{x}_{n}^{\top}{\bf b},$ ${\mathbf b}\in{\R}^p, \; i=1,\ldots,n.$
The corresponding rank scores are
$$a_n(i)=\varphi_n\left(\nfrac{i}{n+1}\right)=I[i=n]-\nfrac 1n, \quad i=1,\ldots,n.$$
Hence,
\begin{eqnarray}\label{11d}
{\mathcal D}_n({\bf b})&=&\sum_{i=1}^n(Y_i-{\bf x}_i^{\top}{\bf b})I[R_{ni}(Y_i-{\mathbf x}_{i}^{\top}{\bf b})=n]
-(\bar{Y}_n-\bar{\bf x}_n^{\top}{\bf b}),\\[2mm]
&=&\max_{1\leq i \leq n}\{Y_i-({\bf x}_i-\bar{\bf x}_n)^{\top}{\bf%
b}\}-\bar{Y}_n=%
\left(Y_i-({\bf x}_i-\bar{\bf x}_n)^{\top}{\bf b}\right)_{n:n}-\bar{Y}_n.\nonumber
\end{eqnarray}
By Jaeckel \cite{Jaeckel1972}, the measure ${\cal D}_n({\mathbf b})$ is continuous, convex and piecewise linear function of ${\mathbf b}\in {\R}^{p}.$ 
The estimate $\widetilde{\bbe}_{nR}^+$ is invariant to the shift in location, thus it estimates only $\beta_1,\ldots,\beta_p.$ 
Define the additional intercept component  $\widetilde{\beta}_{n0}^+$ as
\begin{equation}
\label{11e}
\widetilde{\beta}_{n0}^+=\max\{Y_i-\mathbf{x}_{i}^{\top}%
\widetilde{\bbe}_{nR}^+, \ 1\leq i\leq n\}.
\end{equation}
Denote
\begin{equation}\label{11f}
\widetilde{\bbe}_{nR}^{+*}=\left(
\begin{array}{c}
\widetilde{\beta}_{n0}^+\\[1mm]
\widetilde{\bbe}_{nR}^+\\
\end{array}
\right).
\end{equation}
By \cite{JP2005}, $\widetilde{\bbe}_{nR}^{+*}$ coincides with the maximal regression quantile $\widehat{\bbe}_n^*(1)$. The following theorem shows that $B_n(1)$ is identical to the maximal residual of $Y_i$'s from the R-estimator $\widetilde{\bbe}_{nR}^+.$ It further implies that $B_n(1)$ is identical to the weighted mean of components of the vector $\widetilde{\bbe}_{nR}^{+*}$ with coefficients based on the regressors. 
\begin{thm}\label{Theorem2}
Under the conditions of Theorem \ref{Theorem1},
\begin{description}
	\item[{(i)}] 
\begin{eqnarray}\label{61} 
{B}_n(1)&=&\max_{1\leq i\leq n}\{Y_{ni}-(\mathbf x_{ni}-\overline{\mathbf x}_n)^{\top}\widetilde{\bbe}_{nR}^+\}
=\max_{1\leq i\leq n}(Y_{ni}-\widetilde{Y}_{ni}^R),\\
\mbox{ where } && \widetilde{Y}_{ni}^R=(\mathbf x_{ni}-\overline{\mathbf x}_n)^{\top}\widetilde{\bbe}_{nR}^+, \; i=1,\ldots,n.\nonumber
\end{eqnarray}
\item[{(ii)}]
\begin{equation}
\label{equiv}
{B}_n(1)=\widehat{\beta}_{n0}(1)+\overline{\mathbf x}_n^{\top}\widehat{\bbe}(1)=\widetilde{\beta}_{n0}^++\overline{\mathbf x}_n^{\top}\widetilde{\bbe}_{nR}^+. 
\end{equation} 
\end{description}
\end{thm}
\textbf{Proof.} Indeed, by(\ref{11d}),
$\widetilde{\bbe}_{nR}^+$ minimizes the extreme residual
$(Y_i-(\mathbf{x}_{i}-\bar{\mathbf x}_n)^{\top}{\bf b})_{n:n}$ with respect to
${\bf b}\in \mathbb R^p.$ Hence, by (\ref{11e}) and (\ref{11f}),
$$\widetilde{\beta}_{n0}^++\overline{\mathbf x}_n^{\top}\widetilde{\bbe}_{nR}^+
=\max_{1\leq i\leq n}\{Y_{i}-(\mathbf{x}_{i}-\bar{\mathbf x}_n)^{\top}
\widetilde{\bbe}_{nR}^+\}.$$
On the other hand, the extreme regression quantile minimizes $b_0+\frac 1n\sum_{i=1}^n\mathbf x_i^{\top}\mathbf b$ among all $(b_0,\mathbf b^{\top})^{\top}
\in{\mathbb R}^{p+1}$ such that $b_0+\mathbf x_i^{\top}\mathbf b\geq	Y_i, \; i=1,\ldots,n.$ Because of (\ref{11d}) and (\ref{11e}),  if $b_0+\mathbf x_i^{\top}\mathbf b\geq	Y_i$ for $i=1,\ldots,n,$ then 
$$b_0+\overline{\mathbf x}_n^{\top}\mathbf b\geq \max_{1\leq i\leq n}\{Y_i-b_0-\mathbf x_i^{\top}\mathbf b\}+b_0+\overline{\mathbf x}_n^{\top}\mathbf b=
 \widetilde{\beta}_{n0}^+ +\overline{\mathbf x}_n^{\top}\widetilde{\bbe}_{nR}^+$$ 
what confirms (\ref{equiv}) and hence also (\ref{61}). \hfill $\Box$

\begin{cor} 
Specifically, the AERQ satisfies the following inequalities:
\begin{description}
	\item[(i)] $B_n(1)\leq \max_{i\leq n}\left\{Y_i-(\mathbf x_i-\overline{\mathbf x}_n)^{\top}\mathbf b\right\} \quad \forall \mathbf b\in\mathbb R^p.$ 
	\item[(ii)] $B_n(1)< e_{n:n}+\beta_0+\overline{\mathbf x}_n^{\top}\bbe,$ where $e_{n:n}=\max_{i\leq n}\{e_1,\ldots,e_n\}.$
	\item[(iii)] $B_n(1)\leq \max_{i\leq n}\left\{Y_i-(\mathbf x_i-\bar{x}_n)^{\top}\widehat{\bbe}_n\right\}$ for any estimator of $\bbe$. 
\end{description}	
Hence, the AERQ is equal to maximum of perturbed observations $Y_i, \; i=1,\ldots,n$, where the perturbation is minimum possible caused by the regression with regressors $\mathbf x_1,\ldots,\mathbf x_n.$
\end{cor}
Indeed, the inequalities follow from (\ref{11d}) for any $\mathbf b\in\mathbb R^p,$ and particularly for $\mathbf b=\mathbf 0$ and for any estimate $\widehat{\bbe}$. The sharp inequality in \textit{(ii)} follows from Theorem \ref{Theorem1}. Though $B_n(1)$ is less than $\max_{i\leq n}Y_i,$ the perturbations of $\max_{i\leq n}Y_i$ are minimum possible.\\[2mm]
\textbf{Corollary 2} \textit{The expected conditional shortfall (\ref{shortfall}) of the loss $Y$ under $\bbe^*=\widehat{\bbe}_n^*(1)$ can be alternatively written as}
\begin{eqnarray*}
&&\frac 1n\sum_{i=1}^n\left[Y_i-\mathbf x_i^{*\top}\widehat{\bbe}_n^*(1)\right]
=\overline{Y}_n-\max_{i\leq n}\left\{Y_i-(\mathbf x_i-\overline{\mathbf x}_n)^{\top}\widetilde{\bbe}_{nR}^+\right\}\leq 0
\end{eqnarray*}
\section*{Conclusion}
The extreme (maximal) regression quantile in the linear regression model is a value of the regression parameter which minimizes the sum of the positive residuals. The averaged extreme regression quantile is a useful tool in simultaneous study of extremes and regression which indicates how the regression affects the extreme observations. Two equivalent forms of the AERQ indicate that in the presence of regression it is sharply smaller than $\max_{i\leq n}Y_i$, but the perturbations are the smallest possible over possible values of $\bbe.$ There are promising applications in the analysis of the risk or losses.

\section*{Acknowledgements} The author thanks the Editor, Associated Editor and two Referees for their valuable comments, which helped her to understand better the structure and background of AERQ. The author also thanks Jitka Dupa\v{c}ov\'a for interesting discussions.

\end{document}